\documentclass[10pt]{iopart}
\usepackage{graphicx,bm,amssymb,cite,textcomp,csquotes,color,epsfig}
\eqnobysec

\newcommand{\beq}{\begin{equation}}
\newcommand{\eeq}{\end{equation}}
\newcommand{\beqa}{\begin{eqnarray}}
\newcommand{\eeqa}{\end{eqnarray}}

\newcommand{\abs}[1]{{\left\vert#1\right\vert}}
\newcommand{\ad}{\Delta}
\newcommand{\cu}{{\mathbb Q}}
\newcommand{\dd}{{\rm d}}
\newcommand{\eps}{\varepsilon}
\newcommand{\even}{{\mathop{\rm even}}}
\newcommand{\frad}[2]{\displaystyle{\displaystyle#1\over\displaystyle#2}}
\newcommand{\ii}{{\rm i}}
\newcommand\ljk[2]{\left(\frac{#1}{#2}\right)}
\newcommand{\mod}{\mathop{\rm mod}}
\newcommand{\odd}{{\mathop{\rm odd}}}
\newcommand{\pr}{{\mathop{\rm prime}}}
\newcommand{\Si}{{\mathop{\rm Si}}}
\newcommand{\zed}{{\mathbb Z}}
\renewcommand{\Re}{\mathop{\rm Re}}
\renewcommand{\P}{{\mathcal{P}}}
\newcommand{\Q}{{\mathit{\Pi}}}

\begin{document}

\title[Zeta-regularization and natural boundaries]
{Zeta-regularization and natural boundaries:\\Sums and products of integers and primes}

\author{P L Krapivsky$^{1,2}$ and J M Luck$^3$\footnote{Corresponding author}}

\address{$^1$ Department of Physics, Boston University, Boston, MA 02215, USA}

\address{$^2$ Santa Fe Institute, Santa Fe, NM 87501, USA}

\address{$^3$ Universit\'e Paris-Saclay, CNRS, CEA, Institut de Physique Th\'eorique,
91191~Gif-sur-Yvette, France}

\begin{abstract}
Euler regularized the divergent product of all positive integers and found beautiful formulas for regularized
sums of integer powers of positive integers.
These derivations essentially relied on what is now called the zeta-regularization technique,
although analytical continuation had not yet been invented.
This classic method is however not applicable to the product of all primes,
as the prime zeta function has a natural boundary along the imaginary axis.
Muñoz Garc\'ia and P\'erez-Marco overcame this obstacle and evaluated the product of all primes to $4\pi^2$
by finding an appropriately regularized value of the derivative of the prime zeta function at the origin,
lying on the natural boundary.
Their approach heavily relies on the M\"obius function and
could therefore be called a M\"obius regularization.
We extend this approach in two novel directions.
First, we show how to make sense of the sum of all primes.
This regularization requires going a finite distance beyond the natural boundary.
Second, we determine the regularized products of integers and primes
in the nine imaginary quadratic fields where integers have a unique factorization into primes,
and establish a general power-law relationship between products of integers and primes.
Two well-known examples are Gauss and Eisenstein integers.
The interest in this approach goes beyond number theory.
In a variety of physical situations,
the zeta-regularization technique is indeed not applicable because the relevant zeta function has a natural boundary.
\end{abstract}

\eads{\mailto{pkrapivsky@gmail.com},~\mailto{jean-marc.luck@ipht.fr}}

\maketitle

\section{Introduction and summary of results}

Euler~\cite{Euler} derived the following expression for the product of all positive integers
\beq
\label{Euler-prod}
P= 1\cdot 2 \cdot 3\cdot 4\cdot 5\cdot 6\cdot 7\cdots = \sqrt{2\pi},
\eeq
long before regularized infinite products were put on firm ground
within the so-called zeta-regularization technique,
which is based on analytical continuation
(see~\cite{davis,edwards,vara,varabook} for historical accounts).
Nowadays, zeta-regularized infinite products are used in many fields.
In the physics literature, they are often referred to as renormalized products.
A prominent instance consists in the determinants, i.e., products of all eigenvalues,
of various operators arising in geometry, quantum mechanics and quantum field theory
(see e.g.~\cite{Ray,RS,H,V,SABK,QHS,manin,KV,I,AZ}, and~\cite{Elizalde,EOR,kirsten} for overviews).

In view of~(\ref{Euler-prod}),
it seems logical to wonder whether one can make sense of the product~$\Q$ of all primes.
This question has seemingly been first asked by Soul\'{e}~\cite[p.~101]{SABK}.
It is made difficult by the fact that
the prime zeta function~$\P(s)$ has a natural boundary along the imaginary axis.
It is therefore not analytic in a neighborhood of the origin,
which is required for the applicability of the zeta-regularization approach
to be used in order to express $\Q$ in terms of $\P'(0)$.
Muñoz Garc\'ia and P\'erez-Marco~\cite{MP1,MP2} found a way to overcome this problem and evaluated
the product of primes:
\beq
\label{Soule-prod}
\Q = 2 \cdot 3\cdot 5\cdot 7 \cdot 11 \cdot 13\cdots = (2\pi)^2.
\eeq
The power-law relation
\beq
\Q=P^4
\label{qpnat}
\eeq
between the two regularized products, a curiosity at first sight, actually extends to all number fields
for which we succeeded in computing the analogues of (\ref{Euler-prod}) and (\ref{Soule-prod})
(see~(\ref{qpgal})).

The interest in the derivation of~(\ref{Soule-prod}) goes far beyond the realm of number theory.
A variety of situations can indeed be found in the physics literature,
where the zeta-regularization technique had to be adapted or generalized,
because the relevant zeta function either does not have a meromorphic continuation,
or is not analytic in a neighborhood of $s=0$, or has a natural boundary.
The quantum mechanics of one particle in a one or two-dimensional potential
provides explicit examples where the spectral zeta functions may have poles and branch cuts of any
order in the $s$ plane,
at positions depending on model parameters (see e.g.~\cite{CEZ,FPS}).
More generally, natural boundaries have also been met in statistical mechanics~\cite{NB-Wu,NB-Guttmann},
in gauge theories~\cite{NB-Dunne}, in scattering amplitudes~\cite{NB-Mizera}.
It was recently advocated~\cite{NB-Dunne,NB-Gukov,NB-Adams} that resurgence theory
leads to a new form of unique continuation, beyond analytic continuation.

The formula~(\ref{Soule-prod}) suggests seeking other prime analogues of classical divergent products
and sums admitting a zeta-regularization procedure.
Two especially popular zeta-regularized sums, also going back to Euler, are
\beqa
\label{Euler-sum-1}
1+2+3+4+5+\cdots &=& \zeta(-1)= -\frac{1}{12},\\
\label{Euler-sum-3}
1^3+2^3+3^3+4^3+\cdots &=& \zeta(-3)= \frac{1}{120}.
\eeqa
The sum of positive integers (\ref{Euler-sum-1}) is the most famous zeta-regularized value,
endlessly discussed by science lovers astonished by the fact that this sum is finite and negative.
The sum (\ref{Euler-sum-3}) is the first zeta-regularized value experimentally measured.
It appeared in the Casimir formula~\cite{Casimir} for
the force between two conducting plates in vacuum
(see~\cite{AB} for very recent work on the Casimir effect, and~\cite{BMM} for an overview).

The first challenge we address here is to make sense of the sum of all primes, thereby providing the prime analogue
of (\ref{Euler-sum-1}).
This task looks formidable, because it involves $\P(-1)$, and $s=-1$ sits at a finite distance
beyond the natural boundary.
We then extend the products (\ref{Euler-prod}) and (\ref{Soule-prod}) to the regularized products of integers
and primes in the nine imaginary quadratic fields where integers have a unique factorization into primes.
Hopefully, the techniques developed here in a number-theoretic context can be applied
to some of the physically interesting problems mentioned above.

The present paper is organized as follows.
In section~\ref{nat}, we consider products of integers and primes.
We recall the modern derivation of~(\ref{Euler-prod}) in section~\ref{natint},
whereas in section~\ref{natpri} we recall the derivation of~(\ref{Soule-prod}) by Muñoz Garc\'ia
and P\'erez-Marco and provide further results on the prime zeta function.
The main novel part on primes resides in section~\ref{natsum},
which contains a detailed heuristic derivation of the regularized sum of primes, yielding
\beqa
\P(-1) &=& 2+3+5+7+11+13+\cdots
\nonumber\\
&=&3-\frac{\Si(2\pi)}{\pi}+\sum_{n\ge1}\frac{\mu(n)}{n}\log\zeta(n+1)
\nonumber\\
&=&2.925\,292\,456\,711\,946\dots,
\label{pmun}
\eeqa
where $\mu(n)$ is the M\"obius function, $\zeta(s)$ the Riemann zeta function, and $\Si(x)$ the sine integral function.

The remainder of the paper is devoted to the regularization of the products of integers
and primes in some algebraic extensions of the rationals, namely imaginary quadratic fields where integers
have a unique factorization into primes.
Gauss conjectured that there are exactly nine such quadratic fields with class number one
(see~\cite{Gauss-1,Gauss-high,HW,cohen1,cohen2} for information about this and other class number
problems going back to Gauss).
Sections~\ref{gauss} and~\ref{eisen} are devoted to detailed heuristic computations \`a la Euler,
inspired by~\cite{MP1,MP2}, of the regularized infinite products of Gauss integers $\zed[\ii]$ and Eisenstein
integers $\zed[\omega]$, with $\omega=(-1+\ii\sqrt3)/2=\e^{2\pi\ii/3}$, and of the associated primes.
These two well-known cases are, in several regards,
the most obvious generalizations of the integers and primes.
They are related to the two simplest planar quantum billiards whose spectra are
integrable, namely the square and the equilateral triangle~\cite{IL,I2}.
They also play a part in recent developments in cosmology (see~\cite{CHY} and references therein).

For Gauss integers (see section~\ref{gauss}),
we find that the product $P_4$ of all non-zero Gauss integers and the product $\Q_4$ of all Gauss
primes obeys
\beq
\abs{\Q_4}=\abs{P_4}^8.
\label{qp4}
\eeq
For Eisenstein integers (see section~\ref{eisen}), we obtain similarly
\beq
\abs{\Q_3}=\abs{P_3}^{12}.
\label{qp3}
\eeq
In section~\ref{sib} we extend the computation of the products of integers and primes
to the remaining seven imaginary quadratic fields whose rings of integers have unique factorization.
The absolute discriminants of these fields are the Heegner numbers $\ad=7$, 8, 11, 19, 43, 67, and
163.
There, we find
\beq
\abs{\Q_\ad}=\abs{P_\ad}^4.
\label{qpp}
\eeq

All the above results concerning quadratic number fields can be summarized
into the following universal formula
\beq
\abs{\Q_\ad}=\abs{P_\ad}^{2w_\ad},
\label{qpgal}
\eeq
where the exponent $2w_\ad$ is twice the number of units,
namely 12 for Eisenstein integers (see~(\ref{qp3})),
8 for Gauss integers (see~(\ref{qp4})),
and 4 for the other imaginary quadratic fields (see~(\ref{qpp})),
just as for positive integers and primes (see~(\ref{qpnat})).
Whether the formula~(\ref{qpgal}) holds for more general algebraic number fields remains an open question.
The most obvious extension is to real quadratic fields $\cu\big(\sqrt{d}\big)$,
where $d$ is a positive squarefree integer.
The conjecture by Gauss that there are infinitely many such
real quadratic fields with class number one is still open.

Besides the present number-theoretic context,
a great deal of regularized products of integers,
such as those associated with periodic, aperiodic, and random sequences,
are the subject of current work that will be reported elsewhere~\cite{usinprep}.

\section{Products of integers and primes}
\label{nat}

\subsection{Product of integers}
\label{natint}

The product of all positive integers,
\beq
P=\prod_{n\ge1}n,
\label{pnatural}
\eeq
can be regularized by relating it to the Riemann zeta function
\beq
\zeta(s)=\sum_{n\ge1}n^{-s}=\prod_{p\;\pr}(1-p^{-s})^{-1}.
\label{zeta}
\eeq
We have indeed formally
\beq
P=\exp\sum_{n\ge1}\log n=\exp\left(-\zeta'(0)\right).
\label{plog}
\eeq
This expression is the gist of the zeta-regularization,
consisting in regularizing divergent products such as~(\ref{pnatural})
by means of analytic continuation.
The series and product representations~(\ref{zeta}) of the Riemann zeta function converge for $\Re s >1$
and can be analytically continued to a meromorphic function in the whole complex $s$-plane,
with a simple pole with unit residue at $s=1$, so that
\beq
\lim_{s\to1}(s-1)\zeta(s)=1.
\label{zetalim}
\eeq
In particular, $\zeta(s)$ is analytic in a neighborhood of $s=0$, with
\beq
\zeta(0)=-\frac12,\qquad\zeta'(0)=-\frac12\log 2\pi.
\label{z0}
\eeq
Inserting the second of these expressions into~(\ref{plog}),
we recover Euler's expression~(\ref{Euler-prod}).

\subsection{Product of primes}
\label{natpri}

In the same vein, the product of all primes,
\beq
\Q=\prod_{p\;\pr}p,
\label{qnatural}
\eeq
is formally given by
\beq
\Q=\exp\left(-\P'(0)\right),
\label{qlog}
\eeq
where $\P(s)$ is the prime zeta function (see e.g.~\cite{glaisher,LW,D,froberg}):
\beq
\P(s)=\sum_{p\;\pr}p^{-s}.
\label{zdef}
\eeq

\begin{figure}[!htbp]
\begin{center}
\includegraphics[angle=0,width=.6\linewidth,clip=true]{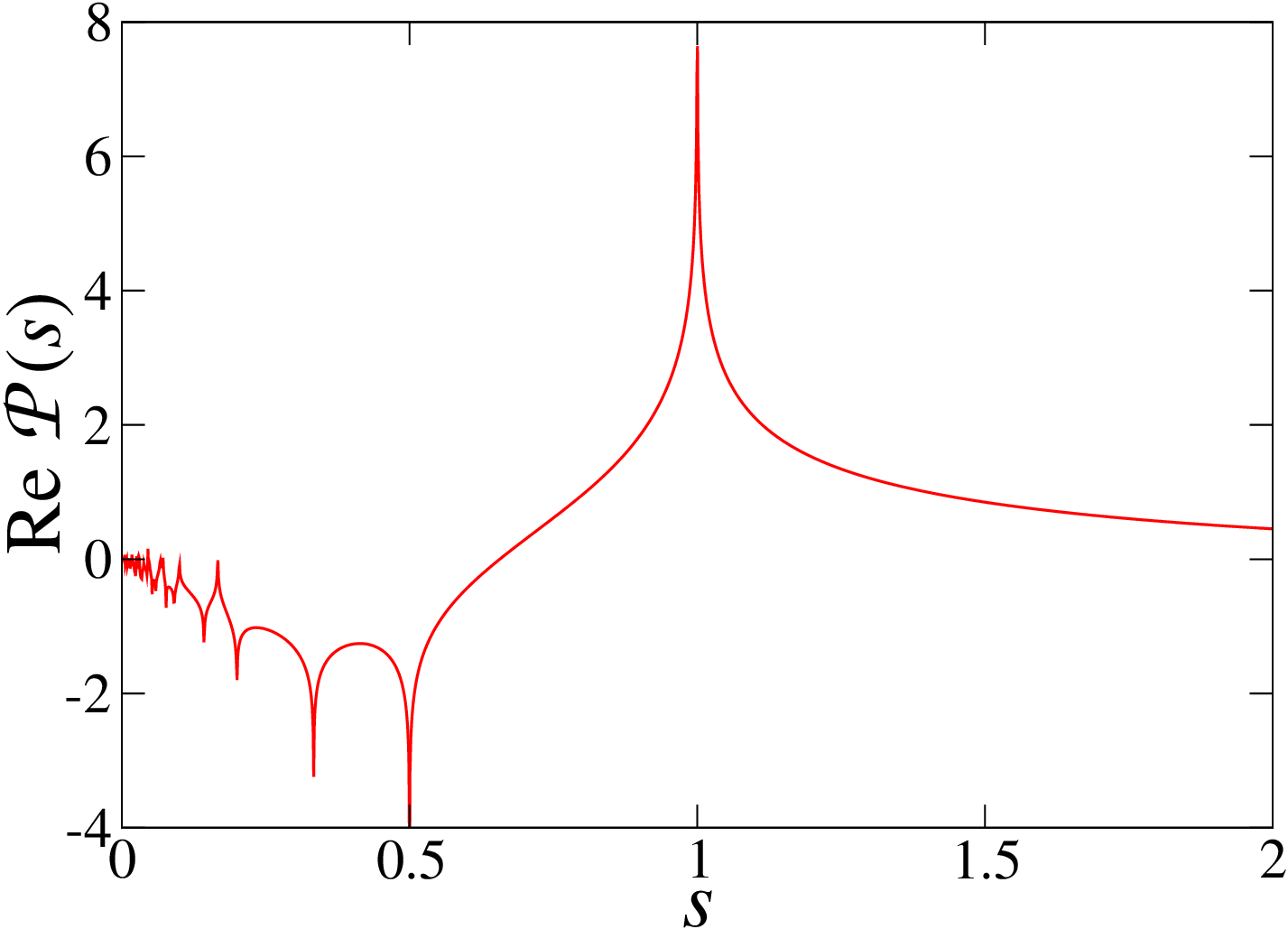}
\caption{Plot of the real part of the prime zeta function for real $s$,
illustrating its infinite sequence of singularities accumulating at the origin (see text).}
\label{primezeta}
\end{center}
\end{figure}

As recalled above, at variance with~(\ref{plog}),
the usage of~(\ref{qlog}) is compromised by the fact that $\P(s)$
has a natural boundary along the imaginary axis
and is therefore not analytic in a neighborhood of $s=0$.
This obstacle has been circumvented by Muñoz Garc\'ia and P\'erez-Marco~\cite{MP1,MP2}.
Their heuristic derivation~\cite{MP1} starts with the Artin-Hasse identity for the exponential
function,
\beq
\e^x=\prod_{n\ge1}(1-x^n)^{-\mu(n)/n},
\label{ah}
\eeq
where $\mu(n)$ is the M\"obius function:
\beq
\mu(n)=\left\{
\begin{array}{cl}
(-1)^k & \hbox{if $n$ is the product of $k$ distinct primes},\\
0 &\hbox{else}.
\end{array}
\right.%}
\label{mudef}
\eeq
Combining~(\ref{zdef}) and~(\ref{ah}), we obtain
\beqa
\exp(\P(s))
&=&\prod_{p\;\pr}\exp(p^{-s})
\nonumber\\
&=&\prod_{p\;\pr}\prod_{n\ge1}(1-p^{-ns})^{-\mu(n)/n}
\nonumber\\
&=&\prod_{n\ge1}\prod_{p\;\pr}(1-p^{-ns})^{-\mu(n)/n}
\nonumber\\
&=&\prod_{n\ge1}\zeta(ns)^{\mu(n)/n}\,.
\label{ezprod}
\eeqa
Taking logarithms of both sides yields
\beq
\P(s)=\sum_{n\ge1}\frac{\mu(n)}{n}\log\zeta(ns).
\label{pmobius}
\eeq
This expression was first established~\cite{glaisher,LW,D,froberg}
by applying the M\"obius inversion formula to the identity
\beq
\log\zeta(s)=\sum_{n\ge1}\frac{\P(ns)}{n},
\eeq
that can be derived from the Euler product formula~(\ref{zeta}).
The expression~(\ref{pmobius}) shows that the prime zeta function has a natural boundary along the imaginary axis.
The non-trivial complex zeros $s_n$ of the Riemann zeta function, belonging to the critical strip ($0<\Re s_n<1$),
indeed yield logarithmic branch cuts accumulating along the whole imaginary axis.
The function $\P(s)$ also has an infinite sequence of logarithmic branch cuts along the real axis
accumulating at the origin (see figure~\ref{primezeta}).
These singularities are located at $s=1/k$ whenever $\mu(k)\ne0$, i.e., $k$ is squarefree.
More specifically, there are upward singularities
at $s=1/k$ whenever $\mu(k)=1$, i.e., $k=1$, 6, 10, 14, 15...
and downward singularities
at $s=1/k$ whenever $\mu(k)=-1$, i.e., $k=2$, 3, 5, 7, 11...

The expression~(\ref{pmobius}) can be used to evaluate regularized values of $\P$ and its
derivatives at $s=0$.
In this approach, sums involving the M\"obius function are consistently evaluated by means
of the zeta-regularization scheme, using
\beq
\label{mu-z}
M(s)=\sum_{n\ge1}\mu(n)\,n^{-s}=\frac{1}{\zeta(s)}\,.
\eeq
In the simplest case of $\P(0)$,~(\ref{pmobius}) yields
\beq
\P(0)=\log\zeta(0)\sum_{n\ge1}\frac{\mu(n)}{n},
\eeq
where the sum over $n$ evaluates to $M(1)=1/\zeta(1)=0$, so that
\beq
\P(0)=0.
\eeq
The value of $\P'(0)$ is of special interest as it enters~(\ref{qlog}).
We have
\beq
\label{Ps-1}
\P'(s)=\sum_{n\ge1}\mu(n)\frac{\zeta'(ns)}{\zeta(ns)},
\eeq
and in particular
\beq
\P'(0)=\frac{\zeta'(0)}{\zeta(0)}\sum_{n\ge1}\mu(n).
\label{zpz}
\eeq
The sum evaluates to $M(0)=1/\zeta(0)$.
We thus obtain, using~(\ref{z0}),
\beq
\P'(0)=\frac{\zeta'(0)}{\zeta(0)^2}=-2\log(2\pi),
\label{zp}
\eeq
and so, using~(\ref{qlog}), we arrive at (\ref{Soule-prod}), the central result of~\cite{MP1,MP2}.

The same approach applies to higher-order derivatives as well.
In the case of $\P''(0)$,
taking the derivative of~(\ref{Ps-1}) and specializing to $s=0$, we arrive at
\beq
\P''(0)=\left(\frac{\zeta''(0)}{\zeta(0)}-\left[\frac{\zeta'(0)}{\zeta(0)}\right]^2\right)
\sum_{n\ge1}n\mu(n).
\eeq
The sum evaluates to $M(-1)=1/\zeta(-1)=-12$.
Using the expression
\beq
\zeta''(0)=\frac{\gamma^2}{2}+\gamma_1-\frac{\pi^2}{24}-\frac{[\log(2\pi)]^2}{2},
\label{z0-2}
\eeq
following from the reflection formula~(\ref{zetarefl}),
we obtain
\beq
\P''(0)=12\gamma^2+24\gamma_1-\pi^2.
\label{P0-2}
\eeq
The constants
\beq
\gamma=\gamma_0=0.577\,215\,664\dots, \qquad \gamma_1=-0.072\,815\,845\dots
\eeq
are the Euler and Stieltjes constants appearing in the Laurent series expansion of the Riemann zeta
function
\beq
\zeta(1+s)=\frac{1}{s}+ \sum_{n\geq 0} \frac{(-1)^n}{n!}\, \gamma_n s^n.
\label{z-exp}
\eeq
A similar computation yields $\P'''(0)=\infty$.
This quantity is indeed proportional to $M(2)=1/\zeta(-2)$, and $\zeta(-2)=0$.

Although evaluating the regularized product $\Q$ of all primes proved to be a challenging problem
resolved only recently,
some infinite products involving primes are easy to regularize.
For instance, the Euler product formula (\ref{zeta}) gives
\beqa
\prod_{p\;\pr}\big(1-p^{2k}\big)=\frac{1}{\zeta(-2k)}=\infty,
\nonumber\\
\prod_{p\;\pr}\big(1-p^{2k-1}\big)=\frac{1}{\zeta(1-2k)}=-\frac{2k}{B_{2k}}
\eeqa
for any integer $k\ge1$,
where $B_{2k}$ are the Bernoulli numbers (see~(\ref{zetaint})).
Specializing these formulas to $k=1$ we obtain
\beq
\prod_{p\;\pr}\big(1-p^{2}\big)=\infty, \quad \prod_{p\;\pr}(1-p\big)=-12.
\eeq
This yields formally
\beq
\prod_{p\;\pr}(1+p)=\infty,
\eeq
and suggests that
\beq
\prod_{p\;\pr}(p+x)
\eeq
has no simple regularized expression.
Such an expression would be a prime analogue of the Lerch formula~\cite{lerch,KW} (see~(\ref{z0x}))
\beq
\label{Lerch:prod}
\prod_{n\ge0}(n+x)=\frac{\sqrt{2\pi}}{\Gamma(x)},
\eeq
which is itself a renormalized avatar of the Weierstrass product formula
\beq
\frac{1}{\Gamma(x)}=x\,\e^{\gamma x}\prod_{n\ge1}(1+x/n)\e^{-x/n}.
\eeq

Finally, we mention that the product
\beq
\prod_{n\ge1}n^n=\exp(-\zeta'(-1))=1.179\,889\,917\dots
\label{nn}
\eeq
admits a prime analogue,
\beq
\prod_{p\;\pr}p^p=\exp(-\P'(-1)),
\label{pp}
\eeq
where
\beq
\label{P-1}
\P'(-1)=\sum_{n\ge1}\mu(n)\frac{\zeta'(-n)}{\zeta(-n)}.
\eeq
The regularized value of this quantity
could be determined along the lines of section~\ref{natsum}.

\section{Sums of integers and primes}
\label{natsum}

\begin{displayquote}
\footnotesize{Why add prime numbers? Prime numbers are made to be multiplied, not added.}

\textit{Lev Landau}
\end{displayquote}

Regularized sums of powers of positive integers are more popular than regularized products.
The well-known examples like (\ref{Euler-sum-1}) and (\ref{Euler-sum-3}) are special cases of the
general Euler's formula~\cite{kks1}
\beq
\label{z-k}
\sum_{n\ge1}n^m=\zeta(-m),
\eeq
expressing the regularized sums of powers of integers in terms of the Bernoulli numbers,
namely
\beq
\zeta(-2k)=0,\quad\zeta(1-2k)=-\frac{B_{2k}}{2k}\qquad(k\ge1).
\label{zetaint}
\eeq

In view of the above, it seems logical to consider the zeta-regularized sums of powers of
primes
\beq
\P(-m)=\sum_{p\;\pr}p^m,
\eeq
for which the formula~(\ref{pmobius}) yields
\beq
\P(-m)=\sum_{n\ge1}\frac{\mu(n)}{n}\log\zeta(-mn).
\label{pm}
\eeq

Investigating the above series is a far more formidable task
than evaluating the regularized product of all primes.
We recall that the prime zeta function has a natural boundary along the imaginary axis.
The regularized product of all primes involves $\P'(0)$,
where the origin sits on the natural boundary,
whereas the regularized sums of powers of primes involve the values of $\P(s)$ at $s=-m$,
at a finite distance beyond the natural boundary.

The famous quote by the physicist Lev Landau, reproduced above,
was made in connection to the Goldbach conjecture
asserting that every even positive integer is the sum of two prime numbers.
Regularizing divergent sums is admittedly a different story.
Finding an application, either in Physics or elsewhere, of $\P(-1)$, or of another $\P(-m)$, would be amusing.
One should keep in mind that more than two hundred years separate
the determination of $\zeta(-3)$ by Euler (see~(\ref{Euler-sum-3}))
and the measurement of the Casimir effect, involving $\zeta(-3)$.

The main goal of this section is to explain
how to make sense of the formula~(\ref{pm}) in the first case of interest,
namely the sum of all primes,
\beq
\P(-1)=\sum_{p\;\pr}p.
\eeq
The following detailed derivation yields the result~(\ref{punana}), already announced in~(\ref{pmun}).
Setting $m=1$ in~(\ref{pm}),
and taking the real parts of the logarithms,
as we expect the regularized sum of all primes to be real,
we obtain
\beq
\P(-1)=\sum_{n\ge1}\frac{\mu(n)}{n}\log\abs{\zeta(-n)}.
\eeq
The Riemann zeta function vanishes at even negative integers (see~(\ref{zetaint})),
so that all terms of the above sum with even $n$ are divergent.
In order to circumvent this, we rather consider~(see~(\ref{pmobius}))
\beq
\P(-1+\eps)=\sum_{n\ge1}\frac{\mu(n)}{n}\log\abs{\zeta(-n+n\eps)}.
\eeq
and investigate the $\eps\to0$ limit of this expression.
To leading order in $\eps$, we have
\beq
\P(-1+\eps)\approx\sum_{n\;\odd}\frac{\mu(n)}{n}\log\abs{\zeta(-n)}
+\sum_{n\;\even}\frac{\mu(n)}{n}\log\abs{n\eps\zeta'(-n)}.
\label{peps}
\eeq
All terms with even $n$ are logarithmically divergent,
whereas terms with odd $n$ are finite.

We therefore have to deal separately with odd and even values of $n$.
To do so, we introduce the partial Dirichlet series
\beq
M_\odd(s)=\sum_{n\;\odd}\mu(n)\,n^{-s},\quad
M_\even(s)=\sum_{n\;\even}\mu(n)\,n^{-s},
\eeq
which can be expressed in terms of the full $M(s)$ introduced in~(\ref{mu-z}).
First, we have the sum rule
\beq
M_\odd(s)+M_\even(s)=M(s)=\frac{1}{\zeta(s)}.
\label{msum}
\eeq
Second, the definition~(\ref{mudef}) of the M\"obius function implies
\beq
\mu(2n)=\left\{
\begin{array}{cl}
-\mu(n) & \hbox{($n$ odd)},\\
0 &\hbox{($n$ even)},
\end{array}\right.%}
\eeq
so that
\beq
M_\even(s)=-\frac{M_\odd(s)}{2^s}.
\label{mrel}
\eeq
Equations~(\ref{msum}) and~(\ref{mrel}) yield
\beq
M_\odd(s)=\frac{2^s}{(2^s-1)\zeta(s)},\quad
M_\even(s)=-\frac{1}{(2^s-1)\zeta(s)}.
\label{mome}
\eeq

The terms proportional to $\log\abs{\eps}$ in~(\ref{peps}) sum up to $M_\even(1)=0$.
This vanishing result is due to the fact that $\zeta(s)$ has a pole at $s=1$ (see~(\ref{zetalim})).
It ensures that~(\ref{peps}) has a finite limit as $\eps\to0$, namely
\beq
\P(-1)=\sum_{n\;\odd}\frac{\mu(n)}{n}\log\abs{\zeta(-n)}
+\sum_{n\;\even}\frac{\mu(n)}{n}\log\abs{n\zeta'(-n)}.
\label{pmun1}
\eeq
The reflection formula for the Riemann zeta function
\beq
\zeta(1-s)=2(2\pi)^{-s}\cos\frac{\pi s}{2}\Gamma(s)\zeta(s)
\label{zetarefl}
\eeq
yields
\beq
\begin{array}{ll}
\abs{\zeta(-n)}=\frad{n!\,\zeta(n+1)}{\pi(2\pi)^n} & \hbox{($n$ odd)},\\
\abs{n\zeta'(-n)}=\frad{n\,n!\,\zeta(n+1)}{2(2\pi)^n} \quad & \hbox{($n$ even)},
\end{array}
\eeq
and so
\beqa
\P(-1)&=&\sum_{n\;\odd}\frac{\mu(n)}{n}\log\frac{n!\,\zeta(n+1)}{\pi(2\pi)^n}
\nonumber\\
&+&\sum_{n\;\even}\frac{\mu(n)}{n}\log\frac{n\,n!\,\zeta(n+1)}{2(2\pi)^n}.
\label{pmun2}
\eeqa
The exact regularization of the above expression goes as follows.
Both series of terms share the common leading asymptotic behavior $(n/(2\pi\e))^n$.
The corresponding regularized sum reads
\beq
\sum_{n\ge1}\mu(n)(\log n-\log(2\pi)-1)=2,
\eeq
where we have used
\beq
\sum_{n\ge1}\mu(n)=\frac{1}{\zeta(0)}=-2,\quad
\sum_{n\ge1}\mu(n)\log n=\frac{\zeta'(0)}{\zeta(0)^2}=-2\log(2\pi).
\eeq
At this stage, we have therefore
\beqa
\P(-1)=2&+&\sum_{n\;\odd}\frac{\mu(n)}{n}\log\frac{n!(\e/n)^n\zeta(n+1)}{\pi}
\nonumber\\
&+&\sum_{n\;\even}\frac{\mu(n)}{n}\log\frac{n\,n!(\e/n)^n\zeta(n+1)}{2}.
\eeqa
The above sums are still divergent.
The arguments of the logarithms of both series of terms indeed respectively behave at large $n$ as
$(2n/\pi)^{1/2}$ and $(\pi n^3/2)^{1/2}$.
The corresponding regularized sums read
\beq
\sum_{n\;\odd}\frac{\mu(n)}{n}\,\frac{\log(2n/\pi)}{2}=-1,\quad
\sum_{n\;\even}\frac{\mu(n)}{n}\,\frac{\log(\pi n^3/2)}{2}=\frac{3}{2},
\eeq
where we have used
\beqa
&&\sum_{n\;\odd}\frac{\mu(n)}{n}=M_\odd(1)=0,\quad
\sum_{n\;\odd}\frac{\mu(n)}{n}\,\log n=-M_\odd'(1)=-2,
\nonumber\\
&&\sum_{n\;\even}\frac{\mu(n)}{n}=M_\even(1)=0,\quad
\sum_{n\;\even}\frac{\mu(n)}{n}\,\log n=-M_\even'(1)=1.
\nonumber\\
\eeqa
We thus obtain the expression
\beq
\P(-1)=\frac{5}{2}+\sum_{n\ge1}\frac{\mu(n)}{n}\log\frac{n!(\e/n)^n}{\sqrt{2\pi n}}
+\sum_{n\ge1}\frac{\mu(n)}{n}\log\zeta(n+1),
\label{pudble}
\eeq
The non-trivial part of the computation involving the regularization of divergent sums
yields the simple rational number $5/2$,
whereas both series are absolutely convergent.
The first of them is evaluated in~\ref{apps}.
Using~(\ref{sires}), we are left with the result
\beq
\P(-1)=3-\frac{\Si(2\pi)}{\pi}+\sum_{n\ge1}\frac{\mu(n)}{n}\log\zeta(n+1),
\label{punana}
\eeq
where $\Si$ denotes the sine integral function.
The series entering~(\ref{punana}) converges exponentially fast, as $\log\zeta(n+1)\approx1/2^{n+1}$.
Summing 60 terms is sufficient to obtain the sum of the series to 15 significant figures,
namely $0.376\,704\,123\,502\,087\dots$
Using~(\ref{sinu}), we obtain
\beq
\P(-1)=2.925\,292\,456\,711\,946\dots
\label{punnum}
\eeq
The results~(\ref{punana}),~(\ref{punnum}) were announced in~(\ref{pmun}).

We emphasize again that $\P(s)$ has a natural boundary at $\Re s=0$,
and so $\P(-1)$ is not the analytic continuation of the prime zeta function,
but the regularized value obtained by means of the identity~(\ref{pmobius}) involving the M\"obius function.
The above procedure could therefore be called a M\"obius regularization.

To close, we mention that the regularized values of $\P'(-1)$ (see~(\ref{P-1})),
and of $\P(-m)$ and $P'(-m)$ for higher integers $m$,
can be evaluated along the very lines of the above derivation.

\section{Gauss integers and primes}
\label{gauss}

Gauss integers form the ring $\zed[\ii]$.
They can be viewed as the vertices of the square lattice,
of the form $z=a+b\ii$, where $a$ and $b$ are integers, and so $\abs{z}^2=a^2+b^2$.

\subsection{Product of Gauss integers}

The product of all Gauss integers is
\beq
P_4=\prod_{(a,b)\ne(0,0)}(a+b\ii).
\label{p4def}
\eeq
The subscript $4$ refers to the absolute discriminant $\Delta=4$.

The ring $\zed[\ii]$ has 4 units, namely $\pm 1$ and $\pm \ii$.
It is invariant under multiplication by $\ii$.
Therefore, assigning a meaning to the phase of $P_4$ is problematic,
and we limit ourselves to
\beq
\abs{P_4}^2=\prod_{(a,b)\ne(0,0)}(a^2+b^2).
\label{p42def}
\eeq
The corresponding Dedekind zeta function,
\beq
\zeta_4(s)=\sum_{(a,b)\ne(0,0)}(a^2+b^2)^{-s},
\label{z42def}
\eeq
is such that
\beq
\abs{P_4}^2=\exp\left(-\zeta_4'(0)\right).
\label{p42}
\eeq
The theory of Dedekind zeta functions is exposed in detail in~\cite{kks2,kks3,cartier} for the
cases of $\zed[\ii]$ and $\zed[\omega]$,
considered in sections~\ref{gauss} and~\ref{eisen}.
In the present situation, we have
\beq
\zeta_4(s)=4\zeta(s)L_4(s),
\label{z4factor}
\eeq
where $\zeta(s)$ is the Riemann zeta function,
and
\beq
L_4(s)=\sum_{n\ge1}\chi_4(n)\,n^{-s}
\eeq
is the Dirichlet $L$-function associated with the Dirichlet character $\chi_4$,
i..e, the multiplicative function of $n \mod 4$ such that
\beq
\chi_4(0)=0,\quad
\chi_4(1)=1,\quad
\chi_4(2)=0,\quad
\chi_4(3)=-1.
\label{chi4def}
\eeq
This $L$-function reads explicitly
\beqa
L_4(s)
&=&\sum_{m\ge0}\left((4m+1)^{-s}-(4m+3)^{-s}\right)
\nonumber\\
&=&\frac{\zeta(s,1/4)-\zeta(s,3/4)}{4^s},
\label{l4h}
\eeqa
in terms of the Hurwitz zeta function
\beq
\zeta(s,x)=\sum_{m\ge0}(m+x)^{-s}.
\label{hurwitz}
\eeq
The analytic structure of the latter function is similar to that of the Riemann zeta function,
which is recovered as $\zeta(s)=\zeta(s,1)$.
In particular~\cite{kks3}
\beq
\zeta(0,x)=\frac12-x,\qquad\zeta'(0,x)=\log\frac{\Gamma(x)}{\sqrt{2\pi}}\,.
\label{z0x}
\eeq
The second expression, and the ensuing regularized infinite product~(\ref{Lerch:prod}),
are known as the Lerch formula~\cite{lerch,KW}.
We also mention that $L_4(s)$ has an Euler product representation of the form
\beq
L_4(s)=\prod_{p\;\pr\ne2}(1-\chi_4(p)\,p^{-s})^{-1},
\eeq
i.e., explicitly
\beq
L_4(s)=\frac{\zeta_{4,1}(s)\zeta_{4,3}(2s)}{\zeta_{4,3}(s)},
\label{l4zz}
\eeq
in terms of the partial prime zeta functions
\beqa
\zeta_{4,1}(s)&=&\prod_{p\;\pr=1\mod4}(1-p^{-s})^{-1},
\nonumber\\
\zeta_{4,3}(s)&=&\prod_{p\;\pr=3\mod4}(1-p^{-s})^{-1},
\label{z4143def}
\eeqa
which obey
\beq
\zeta_{4,1}(s)\zeta_{4,3}(s)=(1-2^{-s})\zeta(s).
\label{z4143z}
\eeq
Using~(\ref{z4factor}) and~(\ref{l4h}), together with~(\ref{z0}) and~(\ref{z0x}), we obtain
\beqa
&&L_4(0)=\frac12,\qquad
L'_4(0)=\log\frac{\Gamma(1/4)}{2\Gamma(3/4)},
\label{l4res}
\\
&&\zeta_4(0)=-1,\qquad
\zeta'_4(0)=\log\frac{2\Gamma(3/4)^2}{\pi\Gamma(1/4)^2}.
\label{zeta4res}
\eeqa
The value of $\zeta_4(0)$ can be interpreted as counting negatively
the point $(0,0)$ that is excluded from the product~(\ref{p42def}) and the sum~(\ref{z42def}).
Finally,~(\ref{p42}) yields
\beq
\abs{P_4}^2=\frac{\pi\Gamma(1/4)^2}{2\Gamma(3/4)^2}=\frac{\Gamma(1/4)^4}{4\pi}=13.750\,371\,636\dots
\label{p4res}
\eeq

\subsection{Product of Gauss primes}

The product of all Gauss primes,
\beq
\Q_4=\prod_{a+b\ii\;\pr}(a+b\ii),
\label{q4}
\eeq
is also invariant under multiplication by $\ii$,
so we consider again
\beq
\abs{\Q_4}^2=\prod_{a+b\ii\;\pr}(a^2+b^2).
\label{q42def}
\eeq
To evaluate this quantity we employ the prime zeta function
\beq
\P_4(s)=\sum_{a+b\ii\;\pr}(a^2+b^2)^{-s},
\eeq
such that
\beq
\abs{\Q_4}^2=\exp\left(-\P_4'(0)\right).
\eeq
Accordingly, in~(\ref{ezprod}) to~(\ref{zpz}), the Riemann zeta function $\zeta(s)$ is to be
replaced by
\beq
Z_4(s)=\prod_{a+b\ii\;\pr}(1-(a^2+b^2)^{-s})^{-1}.
\eeq
It is instructive to determine this function by elementary means
from the sole knowledge of Gauss primes, dating back to Gauss himself.
The Gauss integer $z=a+b\ii$ is prime in the following three situations (see e.g.~\cite[Ch.~XV]{HW}):

\begin{itemize}

\item $z$ is the product of a prime $p=3 \mod 4$ by a unit.
The prime $p$ is said to be inert in $\zed[\ii]$.
It corresponds to 4 distinct Gauss primes, as there are 4 units.

\item $\abs{z}^2=a^2+b^2$ is a prime $p=1 \mod 4$.
The prime $p$ is said to be split.
It corresponds to 8 distinct Gauss primes,
namely the products of $a+\ii b$ and $b+\ii a$ by units.

\item $z$ is the product of $1+\ii$ by a unit.
The prime $\abs{z}^2=2$ is said to be ramified.
It corresponds to 4 distinct Gauss primes.

\end{itemize}

Taking the above multiplicities into account, and using~(\ref{z4143def}), we obtain
\beq
Z_4(s)=\zeta_{4,3}(2s)^4\zeta_{4,1}(s)^8(1-2^{-s})^{-4}.
\label{z4full}
\eeq
Using~(\ref{l4zz}) and~(\ref{z4143z}), this boils down to
\beq
Z_4(s)=\zeta(s)^4L_4(s)^4=\left(\frac{\zeta_4(s)}{4}\right)^4.
\eeq
The analogue of~(\ref{zp}) therefore reads
\beq
\P_4'(0)
=\frac{1}{\zeta(0)}\,\frac{Z_4'(0)}{Z_4(0)}
=\frac{4\zeta_4'(0)}{\zeta(0)\zeta_4(0)}
=8\zeta_4'(0).
\eeq
We thus obtain the result
\beq
\abs{\Q_4}=\abs{P_4}^8,
\eeq
announced in~(\ref{qp4}).

\section{Eisenstein integers and primes}
\label{eisen}

Eisenstein integers form the ring $\zed[\omega]$, with $\omega=(-1+\ii\sqrt3)/2=\e^{2\pi\ii/3}$.
They can be viewed as the vertices of the triangular lattice,
of the form $z=a+b\omega$, where $a$ and $b$ are integers, so that $\abs{z}^2=a^2-ab+b^2$.
The following analysis parallels that of Gauss integers and primes, performed in
section~\ref{gauss}.

\subsection{Product of Eisenstein integers}

The product of all Eisenstein integers is
\beq
P_3=\prod_{(a,b)\ne(0,0)}(a+b\omega).
\label{p3def}
\eeq
The subscript $3$ refers to the absolute discriminant $\Delta=3$.

The ring $\zed[\omega]$ has 6 units, namely $\pm 1$, $\pm\omega$, and $\pm(1+\omega)$.
We again consider
\beq
\abs{P_3}^2=\prod_{(a,b)\ne(0,0)}(a^2-ab+b^2).
\label{p32def}
\eeq
The corresponding Dedekind zeta function,
\beq
\zeta_3(s)=\sum_{(a,b)\ne(0,0)}(a^2-ab+b^2)^{-s},
\label{z32def}
\eeq
such that
\beq
\abs{P_3}^2=\exp\left(-\zeta_3'(0)\right),
\label{p32}
\eeq
reads
\beq
\zeta_3(s)=6\zeta(s)L_3(s),
\label{z3factor}
\eeq
where
\beq
L_3(s)=\sum_{n\ge1}\chi_3(n)\,n^{-s}
\eeq
is the Dirichlet $L$-function associated with the Dirichlet character $\chi_3$,
i..e, the multiplicative function of $n \mod 3$ such that
\beq
\chi_3(0)=0,\quad
\chi_3(1)=1,\quad
\chi_3(2)=-1.
\label{chi3def}
\eeq
This $L$-function reads explicitly
\beqa
L_3(s)
&=&\sum_{m\ge0}\left((3m+1)^{-s}-(3m+2)^{-s}\right)
\nonumber\\
&=&\frac{\zeta(s,1/3)-\zeta(s,2/3)}{3^s}.
\label{l3h}
\eeqa
We also mention that $L_3(s)$ has an Euler product representation of the form
\beq
L_3(s)=\prod_{p\;\pr\ne3}(1-\chi_3(p)\,p^{-s})^{-1},
\eeq
i.e., explicitly
\beq
L_3(s)=\frac{\zeta_{3,1}(s)\zeta_{3,2}(2s)}{\zeta_{3,2}(s)},
\label{l3zz}
\eeq
in terms of the partial prime zeta functions
\beqa
\zeta_{3,1}(s)&=&\prod_{p\;\pr=1\mod3}(1-p^{-s})^{-1},
\nonumber\\
\zeta_{3,2}(s)&=&\prod_{p\;\pr=2\mod3}(1-p^{-s})^{-1},
\label{z3132def}
\eeqa
which obey
\beq
\zeta_{3,1}(s)\zeta_{3,2}(s)=(1-3^{-s})\zeta(s).
\label{z3132z}
\eeq
Using~(\ref{z3factor}) and~(\ref{l3h}), together with~(\ref{z0}) and~(\ref{z0x}), we obtain
\beqa
&&L_3(0)=\frac13,\qquad
L'_3(0)=\log\frac{\Gamma(1/3)}{3^{1/3}\Gamma(2/3)},
\label{l3res}
\\
&&\zeta_3(0)=-1,\qquad
\zeta'_3(0)=\log\frac{3\Gamma(2/3)^3}{2\pi\Gamma(1/3)^3}.
\label{zeta3res}
\eeqa
Finally,~(\ref{p32}) yields
\beq
\abs{P_3}^2=\frac{2\pi\Gamma(1/3)^3}{3\Gamma(2/3)^3}=\frac{\sqrt{3}\,\Gamma(1/3)^6}{4\pi^2}=16.217\,
256\,528\dots
\label{p3res}
\eeq

\subsection{Product of Eisenstein primes}

The product of all Eisenstein primes is
\beq
\Q_3=\prod_{a+b\omega\;\pr}(a+b\omega),
\label{q3}
\eeq
and we again consider
\beq
\abs{\Q_3}^2=\prod_{a+b\omega\;\pr}(a^2-ab+b^2).
\label{q32def}
\eeq
Introducing the prime zeta function
\beq
\P_3(s)=\sum_{a+b\omega\;\pr}(a^2-ab+b^2)^{-s},
\eeq
we have
\beq
\abs{\Q_3}^2=\exp\left(-\P_3'(0)\right).
\eeq
Accordingly, in~(\ref{ezprod}) to~(\ref{zpz}), the Riemann zeta function $\zeta(s)$ is to be
replaced by
\beq
Z_3(s)=\prod_{a+b\omega\;\pr}(1-(a^2-ab+b^2)^{-s})^{-1}.
\eeq
It is again instructive to determine this function by elementary means.
The Eisenstein integer $z=a+b\omega$ is prime in the following three situations (see
e.g.~\cite[Ch.~XV]{HW}):

\begin{itemize}

\item $z$ is the product of a prime $p=2 \mod 3$ by a unit.
The prime $p$ is said to be inert in $\zed[\omega]$.
It corresponds to 6 distinct Eisenstein primes, as there are 6 units.

\item $\abs{z}^2=a^2-ab+b^2$ is a prime $p=1 \mod 3$.
The prime $p$ is said to be split.
It corresponds to 12 distinct Eisenstein primes.

\item $z$ is the product of $1-\omega$ by a unit.
The prime $\abs{z}^2=3$ is said to be ramified.
It corresponds to 6 distinct Eisenstein primes.

\end{itemize}

Taking the above multiplicities into account, and using~(\ref{z3132def}), we obtain
\beq
Z_3(s)=\zeta_{3,2}(2s)^6\zeta_{3,1}(s)^{12}(1-3^{-s})^{-6}.
\label{z3full}
\eeq
Using~(\ref{l3zz}) and~(\ref{z3132z}), this boils down to
\beq
Z_3(s)=\zeta(s)^6L_3(s)^6=\left(\frac{\zeta_3(s)}{6}\right)^6.
\eeq
The analogue of~(\ref{zp}) therefore reads
\beq
\P_3'(0)
=\frac{1}{\zeta(0)}\,\frac{Z_3'(0)}{Z_3(0)}
=\frac{6\zeta_3'(0)}{\zeta(0)\zeta_3(0)}
=12\zeta_3'(0).
\eeq
We thus obtain
\beq
\abs{\Q_3}=\abs{P_3}^{12},
\eeq
announced in~(\ref{qp3}).

\section{Integers and primes in imaginary quadratic fields}
\label{sib}

We now consider a generic imaginary quadratic field $\cu\big(\ii\sqrt{d}\big)$, where $d$ is a
positive squarefree integer.
We shall make extensive use of the material exposed in the books by Cohen~\cite{cohen1,cohen2},
especially in section 10 of~\cite{cohen2}; see also~\cite{HW} for a classic overview.

\subsection{Background}

The ring of integers of the field $\cu\big(\ii\sqrt{d}\big)$ is $\zed[\omega]$, where

\begin{itemize}

\item
If $d=1$ or $2 \mod 4$, we have
\beq
\omega=\ii\sqrt{d},
\label{o1}
\eeq
so that $\abs{\omega}^2=d$.
The corresponding discriminant is
\beq
D=-4d.
\eeq

\item
If $d=3 \mod 4$, we have\footnote{For $d=3$,~(\ref{o2}) yields
$\omega=(1+\ii\sqrt3)/2=\e^{\ii\pi/3}$.
In section~\ref{eisen}, we used the classic notation $\omega$ for $(-1+\ii\sqrt3)/2=\e^{2\ii\pi/3}$.
Despite this change in convention, all results given below are consistent with those of
section~\ref{eisen}.}
\beq
\omega=\frac{1+\ii\sqrt{d}}{2},
\label{o2}
\eeq
so that $\abs{\omega}^2=(d+1)/4$ is again an integer.
The corresponding discriminant is
\beq
D=-d.
\eeq

\end{itemize}

Henceforth we consider rings $\mathbb{Z}[\omega]$
in which every integer uniquely factorizes into a product of primes, up to units and ordering.
In other words, $\zed[\omega]$ must be a unique factorization domain
(or equivalently a principal ideal domain, in the present quadratic realm).
It was conjectured by Gauss, and finally proven by Heegner, Baker and Stark (see~\cite{Gauss-1}),
that only nine imaginary quadratic fields obey this property
(see table~\ref{nine}).
The corresponding absolute discriminants,
\beq
\ad=\abs{D},
\label{dad}
\eeq
are dubbed Heegner numbers.
For the first three of these imaginary quadratic fields, the corresponding rings
$\zed\!\left[\big(1+\ii\sqrt{3}\big)/2\right], \zed[\ii]$, and
$\zed\!\left[\big(1+\ii\sqrt{7}\big)/2\right]$
are Eisenstein integers, Gauss integers, and Klein integers.

The number of units (i.e., of roots of unity) in $\zed[\omega]$ is denoted by $w_\ad$.
The first smaller values of $\ad$ correspond to Euclidean domains, whereas the last four do not.
The first two cases respectively correspond to Eisenstein and Gauss integers,
investigated in detail in section~\ref{eisen} ($\ad=3$) and section~\ref{gauss} ($\ad=4$).
These are the only examples where $\abs{\omega}^2=1$.
This property allows for extra symmetries (those of the triangular and square lattices),
and for higher numbers of units ($w_3=6$ and $w_4=4$), whereas $w_\ad=2$ in the generic case.
The following developments hold for any of the nine imaginary quadratic fields listed above.
Notations are consistent with those of previous sections.

\begin{table}[!ht]
\begin{center}
$$
\begin{array}{|c||c|c|c|c|c|c|c|c|c|}
\hline
\ad & 3 & 4 & 7 & 8 & 11 & 19 & 43 & 67 & 163 \\
\hline
d & 3 & 1 & 7 & 2 & 11 & 19 & 43 & 67 & 163 \\
\hline
w_\ad & 6 & 4 & 2 & 2 & 2 & 2 & 2 & 2 & 2 \\
\hline
\abs{\omega}^2 & 1 & 1 & 2 & 2 & 3 & 5 & 11 & 17 & 41 \\
%\mu & 0 & 12 & -15 & 20 & -32 & -96 & -960 & -5280 & -640320 \\
\hline
\end{array}
$$
\caption{The nine imaginary quadratic fields $\cu(\ii\sqrt{d})$ with unique factorization: absolute
discriminant $\ad$,
defining integer $d$, number $w_\ad$ of units, and integer value of~$\abs{\omega}^2$.}
\label{nine}
\end{center}
\end{table}

\subsection{Product of integers in $\cu\big(\ii\sqrt{d}\big)$}

The product of all integers in $\cu\big(\ii\sqrt{d}\big)$ is
\beq
P_\ad=\prod_{(a,b)\ne(0,0)}(a+b\omega).
\label{pddef}
\eeq
We again rather consider
\beq
\abs{P_\ad}^2=\prod_{(a,b)\ne(0,0)}\abs{a+b\omega}^2.
\label{pd2def}
\eeq
The corresponding Dedekind zeta function,
\beq
\zeta_\ad(s)=\sum_{(a,b)\ne(0,0)}\abs{a+b\omega}^{-2s},
\label{zd2def}
\eeq
such that
\beq
\abs{P_\ad}^2=\exp\left(-\zeta_\ad'(0)\right),
\label{pd2}
\eeq
reads
\beq
\zeta_\ad(s)=w_\ad\zeta(s)L_\ad(s),
\label{zdfactor}
\eeq
where $w_\ad$ is the number of units, $\zeta(s)$ is the Riemann zeta function, and
\beq
L_\ad(s)=\sum_{n\ge1}\chi_\ad(n)\,n^{-s}
\eeq
is the Dirichlet $L$-function associated with the character
\beq
\chi_\ad(n)=\ljk{D}{n},
\eeq
where $D=-\ad<0$ (see~(\ref{dad})), and $\ljk{D}{n}$ is the Kronecker symbol.
The connection with the derivations of~(\ref{z4full}) and~(\ref{z3full})
is made by recalling that a $p$ is inert, split, or ramified in $Z[\omega]$
according to whether $\ljk{D}{p}$ equals $-1$, 0, or $+1$.

We have
\beqa
L_\ad(0)&=&\frac{2}{w_\ad},
\label{lz}
\\
L'_\ad(0)&=&-\frac{2\log\ad}{w_\ad}+\sum_{n=1}^\ad\chi_\ad(n)\,\log\Gamma\!\left(\frac{n}{\ad}\right
).
\label{lpz}
\eeqa
Recalling~(\ref{z0}) and inserting~(\ref{lz}) into~(\ref{zdfactor}) yields
\beq
\zeta_\ad(0)=-1.
\eeq
This result can again be interpreted as counting negatively
the point $(0,0)$ that is excluded from the product~(\ref{pd2def}) and the sum~(\ref{zd2def}).

Furthermore, inserting~(\ref{zdfactor}) into~(\ref{pd2}) yields
\beq
\abs{P_\ad}^2=\exp\left(\frac{w_\ad}{2}\left(L_\ad(0)\log2\pi+L'_\ad(0)\right)\right).
\eeq
Inserting~(\ref{lz}) and~(\ref{lpz}) into the latter formula, we obtain
\beq
\abs{P_\ad}^2=\frac{2\pi}{\ad}\left(\prod_{n=1}^\ad\Gamma\!\left(\frac{n}{\ad}\right)^{\chi_\ad(n)}\right)^{w_\ad/2}.
\label{cs}
\eeq
The product in the right-hand side is known as the Chowla-Selberg gamma ratio.
The above result can be recast as
\beq
\abs{P_\ad}^2=4\pi^2\abs{\eta(\omega)}^4,
\label{etaiden}
\eeq
where $\eta(\omega)$ denotes Dedekind's eta function, namely
\beq
\eta(\omega)=q^{1/24}\prod_{n\ge1}(1-q^n),\qquad q=\e^{2\pi\ii\omega},
\eeq
so that $q=\e^{-\pi\sqrt{\ad}}$ for $\ad=4$ or 8, and $q=-\e^{-\pi\sqrt{\ad}}$ in the other cases
(see~(\ref{o1}),~(\ref{o2})).

\subsection{Product of primes in $\cu\big(\ii\sqrt{d}\big)$}

The product of all primes in $\cu\big(\ii\sqrt{d}\big)$ is
\beq
\Q_\ad=\prod_{a+b\omega\;\pr}(a+b\omega).
\label{qd}
\eeq
We again rather consider
\beq
\abs{\Q_\ad}^2=\prod_{a+b\omega\;\pr}\abs{a+b\omega}^2.
\label{qd2def}
\eeq
Introducing the prime zeta function
\beq
\P_\ad(s)=\sum_{a+b\omega\;\pr}\abs{a+b\omega}^{-2s},
\eeq
we obtain
\beq
\abs{\Q_\ad}^2=\exp\left(-\P_\ad'(0)\right).
\eeq
Accordingly, in~(\ref{ezprod}) to~(\ref{zpz}), the Riemann zeta function $\zeta(s)$ is to be
replaced by
\beq
Z_\ad(s)=\prod_{a+b\omega\;\pr}(1-\abs{a+b\omega}^{-2s})^{-1}.
\eeq
This function reads (see~(\ref{zdfactor}))
\beq
Z_\ad(s)=\left(\zeta(s)L_\ad(s)\right)^{w_\ad}=\left(\frac{\zeta_\ad(s)}{w_\ad}\right)^{w_\ad}.
\eeq
The analogue of~(\ref{zp}) therefore reads
\beq
\P_\ad'(0)
=\frac{1}{\zeta(0)}\,\frac{Z_\ad'(0)}{Z_\ad(0)}
=\frac{w_\ad\zeta_\ad'(0)}{\zeta(0)\zeta_\ad(0)}
=2w_\ad\zeta_\ad'(0).
\eeq
We thus obtain the general formula
\beq
\abs{\Q_\ad}=\abs{P_\ad}^{2w_\ad},
\eeq
which is~(\ref{qpgal}).
The exponent relating the regularized products $\abs{P_\ad}$ of integers and $\abs{\Q_\ad}$ of
primes
is nothing but twice the number $w_\ad$ of units, i.e.,
12 for Eisenstein integers (see~(\ref{qp3})),
8 for Gauss integers (see~(\ref{qp4})),
and 4 for the other imaginary quadratic fields (see~(\ref{qpp})),
just as for usual integers and primes (see~(\ref{qpnat})).

\subsection{Explicit formulas and numerical values}

Here we give a list of expressions for $\abs{P_\ad}^2$,
as given by the Chowla-Selberg gamma ratios~(\ref{cs}).
Explicit formulas are given for the five Euclidean cases,
and only numerical values for the four non-Euclidean ones,
where explicit formulas are too lengthy to be reported here.
The first two agree with~(\ref{p3res}) and~(\ref{p4res}), as should be.
\beqa
\abs{P_3}^2
&=&
\frac{2\pi}{3}\,\frac{\Gamma(1/3)^3}{\Gamma(2/3)^3}=16.217\,256\,528\dots,
\nonumber\\
\abs{P_4}^2
&=&\frac{\pi}{2}\,\frac{\Gamma(1/4)^2}{\Gamma(3/4)^2}=13.750\,371\,636\dots,
\nonumber\\
\abs{P_7}^2
&=&\frac{2\pi}{7}\,\frac{\Gamma(1/7)\Gamma(2/7)\Gamma(4/7)}{\Gamma(3/7)\Gamma(5/7)\Gamma(6/7)}=9.889
\,009\,200\dots,
\nonumber\\
\abs{P_8}^2
&=&\frac{\pi}{4}\,\frac{\Gamma(1/8)\Gamma(3/8)}{\Gamma(5/8)\Gamma(7/8)}=8.973\,175\,814\dots,
\nonumber\\
\abs{P_{11}}^2
&=&\frac{2\pi}{11}\,\frac{\Gamma(1/11)\Gamma(3/11)\Gamma(4/11)\Gamma(5/11)\Gamma(9/11)}{\Gamma(2/11)
\Gamma(6/11)\Gamma(7/11)\Gamma(8/11)\Gamma(10/11)}
\nonumber\\
&=&6.953\,831\,684\dots,
\nonumber\\
\abs{P_{19}}^2&=&4.028\,703\,050\dots,
\nonumber\\
\abs{P_{43}}^2&=&1.274\,160\,276\dots,
\nonumber\\
\abs{P_{67}}^2&=&0.543\,303\,840\dots,
\nonumber\\
\abs{P_{163}}^2&=&0.049\,335\,689\dots
\eeqa

The above numerical values exhibit a fast decay as a function of the absolute discriminant $\ad$.
This observation can be made quantitative by means of~(\ref{etaiden}), yielding
\beq
\abs{P_\ad}^2\approx4\pi^2\abs{q}^{1/6}\approx4\pi^2\e^{-\pi\sqrt{\ad}/6}
\eeq
when $\ad$ is large, up to tiny corrections proportional to $\abs{q}=\e^{-\pi\sqrt{\ad}}$.

\ack

It is a pleasure to thank Michel Bauer, Henri Cohen and Ricardo P\'erez-Marco for useful exchanges.

\appendix

\section{Evaluation of a series entering~(\ref{pudble})}
\label{apps}

This appendix is devoted to the evaluation of the series
\beq
S=\sum_{n\ge1}\frac{\mu(n)}{n}\log\frac{n!(\e/n)^n}{\sqrt{2\pi n}}
\eeq
which enters~(\ref{pudble}).
The derivation proceeds in four steps.

\begin{itemize}

\item
Using the asymptotic expansion
\beq
\log n!=n\log n-n+\log\sqrt{2\pi n}+\sum_{k\ge1}\frac{B_{2k}}{2k(2k-1)n^{2k-1}},
\eeq
where $B_{2k}$ are again the Bernoulli numbers, as well as~(\ref{pmobius}), we get
\beq
S=\sum_{n\ge1}\frac{\mu(n)}{n}\sum_{k\ge1}\frac{B_{2k}}{2k(2k-1)n^{2k-1}}
=\sum_{k\ge1}\frac{B_{2k}}{2k(2k-1)\zeta(2k)}.
\eeq

\item
Using the identity
\beq
B_{2k}=(-1)^{k-1}\frac{2(2k)!}{(2\pi)^{2k}}\,\zeta(2k),
\eeq
we obtain
\beq
S=2\sum_{k\ge1}(-1)^{k-1}\frac{(2k-2)!}{(2\pi)^{2k}}
=\frac{1}{2\pi^2}\sum_{m\ge0}(-1)^m\frac{(2m)!}{(2\pi)^{2m}}.
\eeq
The series in the rightmost side is divergent but alternating.

\item
Using the integral representation
\beq
(2m)!=\int_0^\infty t^{2m}\e^{-t}\,\dd t,
\eeq
we obtain some kind of a Borel sum of the above series, namely
\beq
S=\frac{1}{2\pi^2}\int_0^\infty\e^{-t}\,\dd t\sum_{m\ge0}\left(-\frac{t^2}{4\pi^2}\right)^m
=2\int_0^\infty\frac{\e^{-t}}{4\pi^2+t^2}\,\dd t.
\eeq

\item
Using~\cite[Eq.~3.354.1]{GR}, the above expression evaluates to
\beq
S=\frac{1}{2}-\frac{\Si(2\pi)}{\pi},
\label{sires}
\eeq
where
\beq
\Si(x)=\int_0^x\frac{\sin u}{u}\,\dd u
\eeq
is the sine integral function.
The very accurate numerical value
\beq
\Si(2\pi)=1.418\,151\,576\,132\,628\dots
\label{sinu}
\eeq
can be found in~\cite{OEIS} as sequence number A101512.

\end{itemize}

\section*{References}

\bibliography{products.bib}

\end{document}